\documentclass{elsart}  
\usepackage{ucs}
\usepackage[utf8x]{inputenc}
\PrerenderUnicode{é}

\usepackage{amssymb}

\usepackage[english,francais]{babel}

\newtheorem{theorem}{Theorem}[section]

\newtheorem{e-proposition}[theorem]{Proposition}

\newtheorem{e-definition}[theorem]{Definition\rm}

\newtheorem{theoreme}{Th\'eor\`eme}[section]

\renewcommand{\theequation}{\arabic{equation}}
\def\og{\leavevmode\raise.3ex\hbox{$\scriptscriptstyle\langle\!\langle$~}}
\def\fg{\leavevmode\raise.3ex\hbox{~$\!\scriptscriptstyle\,\rangle\!\rangle$}}
\def\di{\displaystyle}

\def\R{\mathbb R}

\journal{the Acad\'emie des sciences}
\begin{document}
\centerline{}
\begin{frontmatter}




%
\selectlanguage{francais}
\title{Propagation de fronts dans les \'equations de Fisher-KPP avec diffusion fractionnaire}



\author[a]{Xavier Cabré},
\ead{xavier.cabre@upc.edu}
\author[b]{Jean-Michel Roquejoffre}
\ead{roque@mip.ups-tlse}
\address[a]{ICREA et
Universitat Politècnica de Catalunya,
Dep. de Matemàtica Aplicada I,
Av. Diagonal 647,
08028 Barcelone, Espagne}
\address[b]{Institut de Math\'ematiques, Univ. de Toulouse et CNRS (UMR 5219),
118 route de Narbonne, 31062 Toulouse, France}




\medskip
\selectlanguage{francais}

\begin{abstract}
\selectlanguage{francais}
On s'intéresse dans cette note à l'équation de Fisher-KPP dans l'espace entier,
où le laplacien est remplacé par le générateur d'un semi-groupe de Feller à noyau lentement
décroissant, un exemple important étant le laplacien fractionnaire. A la différence
de l'équation de Fisher-KPP classique, où l'état stable envahit l'état instable à une
vitesse constante en temps, nous montrons que la vitesse d'invasion est exponentielle en temps.
Ces résultats apportent une justification mathématiquement rigoureuse à de nombreuses heuristiques
sur ce modèle.
\vskip 0.5\baselineskip

\selectlanguage{english}
\noindent{\bf Abstract}
\vskip 0.5\baselineskip
\noindent
{\bf Front propagation in Fisher-KPP equations with fractional diffusion. }
We study in this note the Fisher-KPP equation where the Laplacian is replaced by the
generator of a Feller semigroup with slowly decaying kernel, an important example being
the fractional Laplacian. Contrary to what happens in the standard Laplacian case, where the
stable state invades the unstable one at constant speed, we prove here that invasion holds
at an exponential in time velocity. These results provide a mathematically rigorous justification
of numerous heuristics about this model.

\end{abstract}
\end{frontmatter}

\selectlanguage{english}
\section*{Abridged English version}
\noindent
Consider the reaction-diffusion equation
\begin{equation}
\label{e0.1}
u_t+Au=f(u)\ (t>0,\ x\in\R^N),\ \ \ \ \ u(0,x)\in[0,1].
\end{equation}
where $A$ is the generator of a Feller semigroup (see \cite{BCP}). That is, 
the solution of $v_t+Av=0$, $v(0,x)=v_0(x)$ is given by
$$
e^{-tA}v_0=p(t,\cdot )*v_0,
$$
where $p$ is a positive kernel satisfying assumptions (i)-(iv) below. In the whole note we will assume
\begin{equation}
\label{e0.2}
\frac{B^{-1}}{t^{\frac{N}{2\alpha}}(1+\vert t^{-{\frac1{2\alpha}}}x\vert^{N+2\alpha})}\leq p(t,x)\leq
\frac{B}{t^{\frac{N}{2\alpha}}(1+\vert t^{-{\frac1{2\alpha}}}x\vert^{N+2\alpha})},
\qquad\quad 0<\alpha<1,
\end{equation}
where $B$ is a constant.
An important example for $A$
is the fractional Laplacian
$$
(-\Delta)^\alpha u(x)=c_{N,\alpha}\int_{\R^N}\frac{u(x)-u(y)}{\vert x-y\vert^{N+2\alpha}}\ dy,
$$
where the above expression has to be understood in the sense of principal values. 
The function
$f$ is concave, and satisfies $f(0)=f(1)=0$ and $f'(1)<0<f'(0)$. 
When $A=-\Delta$ (i.e. $\alpha =1$), equation (\ref{e0.1}) is known as the 
Fisher-KPP equation \cite{KPP}.

We wish to understand the large-time behaviour of (\ref{e0.1}). In view of the 
classical reference \cite{AW} - which treats the case $\alpha=1$ - we expect the stable state 1 to invade the unstable state 0,
and we wish to know at which rate. In \cite{AW}, it is linear in time. Numerous heuristics from physics papers
- see for instance \cite{MVV} - suggest that, in the case $\alpha\in (0,1)$, 
the invasion rate will be exponential in time.
Our results provide a mathematically rigorous proof of this.

The first result concerns compactly supported initial data.
\begin{theorem} 
\label{t0.1}
Let $\alpha\in(0,1)$ and set $c_*=\di\frac{f'(0)}{N+2\alpha}$.
Assume $u(0,\cdot )\not\equiv 0$ to be compactly supported. Then: 

\noindent - for all
$c<c_*$, we have $\di\lim_{t\to+\infty}\inf_{\vert x\vert\leq e^{ct}}u(t,x)=1$,

\noindent - for all $c>c_*$, we have $\di\lim_{t\to+\infty}\sup_{\vert x\vert\geq e^{ct}}u(t,x)=0$.
\end{theorem}
The following result deals with nondecreasing initial data.
\begin{theorem} 
\label{t0.2}
Let $\alpha\in(0,1)$ and $N=1$. Set $c_{**}=\di\frac{f'(0)}{2\alpha}$.
Assume $u(0,\cdot )\not\equiv 0$ 
to be nondecreasing and such that $\R_-\cap{\rm{supp}}(u(0,\cdot ))$ is compact. Then:

\noindent - for all $c<c_{**}$, we have  $\di\lim_{t\to+\infty}\inf_{x\geq -e^{ct}}u(t,x)=1$,

\noindent - for all
$c>c_{**}$, we have $\di\lim_{t\to+\infty}\sup_{x\leq -e^{ct}}u(t,x)=0$.
\end{theorem}
The proofs of these theorems also give that (\ref{e0.1}) cannot have travelling wave 
solutions if $\alpha<1$, 
in sharp contrast with the case $\alpha=1$.

In the case $\alpha=1/2$ and $N=1$, we have sharper asymptotics on the dynamics of the level sets of solutions with compactly supported initial data as in Theorem~\ref{t0.1}.
\begin{theorem} 
\label{t0.3}
Let $\alpha=\di\frac12$, $N=1$, and $f(u)=u(1-u)$.
Assume $u(0,\cdot )\not\equiv 0$ to be compactly supported. For $\lambda\in(0,1)$, set 
$$
x_\lambda^-(t)=\inf\{x:\ u(t,x)=\lambda\} \ \ \ \hbox{and}\ \ \ x_\lambda^+(t)=\sup\{x:\ u(t,x)=\lambda\}.
$$
Then, there exists a constant $C_\lambda>1$ such that, for $t$ large enough,
$$
\di{-C_\lambda e^{t/2}\leq x_\lambda^-(t)\leq -\frac1{C_\lambda}e^{t/2}}\ \ \ \hbox{and}\ \ \ \di{\frac1{C_\lambda}e^{t/2}\leq x_\lambda^+(t)\leq C_\lambda e^{t/2}}.
$$
\end{theorem}

All these results will be established and developed in \cite{CR}.
\newpage

\selectlanguage{francais}

\section{Introduction, motivation}
\noindent Soit $A$ le générateur d'un semi-groupe de Feller (voir \cite{BCP}), 
i.e. un opérateur sur l'espace des
fonctions continues bornées sur $\R^N$ tel que le problème d'évolution
\begin{equation}
\label{e1.0}
v_t+Av=0\ (t>0,\ x\in\R^N),\ \ \ v(0,x)=v_0(x)
\end{equation}
admette une unique solution $v(t,x)$ pouvant être représentée sous la forme
$$
v(t,x)=\int_{\R^N}p(t,x-y)v_0(y)\ dy,
$$
où le noyau $p$ vérifie les propriétés suivantes: 

(i). $p(t,x)\geq 0$ et $p(t,\cdot )$ est de masse unité, i.e. $\di\int p(t,x)\ dx=1$;

(ii). $p$ régulière sur $\R_+\times\R^N$;

(iii). $p(t,\cdot )*p(s,\cdot )=p(t+s, \cdot )$ pour tout couple $(s,t)\in\R_+^2$;

(iv). il existe $B>1$ tel que, pour tous $(t,x)\in
\R_+\times\R^N$ on ait:
\begin{equation}
\label{e1.2}
\frac{B^{-1}}{t^{\frac{N}{2\alpha}}(1+\vert t^{-{\frac1{2\alpha}}}x\vert^{N+2\alpha})}\leq p(t,x)\leq
\frac{B}{t^{\frac{N}{2\alpha}}(1+\vert t^{-{\frac1{2\alpha}}}x\vert^{N+2\alpha})},
\qquad\quad 0<\alpha<1.
\end{equation}

Notons que si $A$ est le générateur d'un processus de Markov  stable,
i.e. il existe une fonction $a(t)$ telle que $p(t,x)=a(t)^{-N}p(1,x/ a(t))$,
alors il existe $\alpha\in(0,1]$ tel que $a(t)=t^{\frac1{2\alpha}}$ - comme dans 	
l'expression (\ref{e1.2}).
Un exemple important d'un tel opérateur est le laplacien fractionnaire $(-\Delta)^\alpha$, avec 
$0<\alpha<1$: si $u$ est assez régulière et à croissance assez lente à l'infini, 
$$
(-\Delta)^\alpha u(x)=c_{N,\alpha}\int_{\R^N}\frac{u(x)-u(y)}{\vert x-y\vert^{N+2\alpha}}\ dy,
\qquad\quad 0<\alpha<1,
$$
où l'expression ci-dessus est à prendre au sens des valeurs principales,
et où la constante $c_{N,\alpha}$ est ajustée pour que le symbole de $(-\Delta)^\alpha$ soit
$\vert\xi\vert^{2\alpha}$. Voir par exemple \cite{Du}. 
Le cas $\alpha=1/2$ est un cas important où $p(t,x)$ admet une expression explicite:
$$
p(t,x)=\frac{B_Nt}{(t^{2}+\vert x\vert^2)^{(N+1)/2}}
=\frac{B_N}{t^N(1+\vert t^{-1}x\vert^2)^{(N+1)/2}},
$$
où $B_N$ est choisie pour assurer la propriété (i).

Soit $f$ une fonction de classe $C^1$ sur $[0,1]$, concave, telle que $f(0)=f(1)=0$ et
$f'(1)<0<f'(0)$ - par exemple $f(u)=u(1-u)$. On s'intéresse au comportement en grand temps
des solutions $u(t,x)$ du problème d'évolution
\begin{equation}
\label{e1.1}
u_t+Au=f(u)\ (t>0,\ x\in\R^N),\ \ \ \ \ u(0,x)\in[0,1].
\end{equation}
Quand $A=-\Delta$ (i.e. $\alpha=1$), 
(\ref{e1.1}) est connue sous le nom d'équation de Fisher-KPP
\cite{KPP}. 
Le problème de Cauchy pour (\ref{e1.1}) est, de par les propriétés (i)-(iii) de $A$,
trivialement bien posé sur $\R_+\times\R^N$, et vérifie de plus le principe de 
comparaison suivant: pour deux solutions $u$ et $v$ de (\ref{e1.1}),
si $u(0,\cdot )\leq v(0,\cdot )$  alors $u\leq v$ sur $\R_+\times\R^N$. 
On s'attend à ce que l'état stable
1 envahisse l'état instable 0, et on souhaite savoir à quelle vitesse.

Ce type d'équation de réaction-diffusion apparaît dans de nombreux
modèles de physique (turbulence, plasmas chauds, modèles de flammes) lorsque les
phénomènes diffusifs ne sont pas correctement décrits par des processus gaussiens.
Voir par exemple 
\cite{MVV} pour une description de quelques-uns de ces modèles. Bien qu'heuristiques, ceux-ci
peuvent apporter des informations qualitatives. L'équation (\ref{e1.1}) apparaît aussi dans une
classe de modèles de dynamique des populations, et peut s'obtenir, dans un certain régime spatio-temporel,
comme l'asymptotique d'un modèle intégro-différentiel - voir \cite{BRR}. L'équation (\ref{e1.1})
avec $A=-\Delta$ (i.e. $\alpha=1$) peut s'obtenir comme asymptotique du même modèle, dans un régime spatio-temporel
différent - voir \cite{KPP}. Il est donc légitime de chercher à comprendre les propriétés qualitatives
de (\ref{e1.1}). De fait, de nombreuses heuristiques physiques existent, voir là encore \cite{MVV}. Le but
de ce travail est d'apporter une justification rigoureuse sur le plan mathématique.

\section{Résultats principaux}
\noindent Lorsque $\alpha=1$ (équation de Fisher-KPP) on rappelle que l'équation (\ref{e1.1}) devient
\begin{equation}
\label{e2.1}
u_t-\Delta u=f(u)\ (t>0,\ x\in\R^N),\ \ \ \ \ u(0,x)\in[0,1].
\end{equation}
Un résultat très général sur (\ref{e2.1}) est dû à Aronson-Weinberger \cite{AW};
il décrit le devenir d'une perturbation à support compact.
\begin{theoreme}[\cite{AW}]
Soit $u$ solution de {\rm (\ref{e2.1})} et supposons que $u(0,\cdot )\not\equiv 0$ est
à support compact. Considé-\-rons la quantité $c_{*,1}=2\sqrt{f'(0)}$. Alors:

\noindent - pour tout $c<c_{*,1}$, on a $\di\lim_{t\to+\infty}\inf_{\vert x\vert\leq ct}u(t,x)=1$,

\noindent -  pour tout $c>c_{*,1}$, on a $\di\lim_{t\to+\infty}\sup_{\vert x\vert\geq ct}u(t,x)=0$.
\end{theoreme}
Il se trouve que (\ref{e2.1}) admet des solutions d'ondes progressives monodimensionnelles
connectant 0 à 1, i.e. des solutions de la forme $\phi(x\cdot e+ct)$, avec
$$
-\phi''+c\phi'=f(\phi), \ \ \ \phi(-\infty)=0,\ \phi(+\infty)=1. 
$$
Le réel $c_{*,1}$ est justement la plus petite vitesse d'onde possible, et Komogorov, Petrovskii et Piskunov
\cite{KPP} ont démontré que la solution de (\ref{e2.1}) avec $N=1$, issue de la fonction de Heaviside 
$H(x)$, converge en grand temps vers un profil d'onde de vitesse $c_{*,1}$.

Nos résultats sont les suivants.
\begin{theoreme} 
\label{t1.1}
Soit $\alpha\in(0,1)$, et considérons la quantité $c_*=\di\frac{f'(0)}{N+2\alpha}$.
Soit $u$ solution de  {\rm (\ref{e1.1})} et supposons que $u(0,\cdot )\not\equiv 0$ 
est à support compact. Alors:

\noindent - pour tout $c<c_*$, on a $\di\lim_{t\to+\infty}\inf_{\vert x\vert\leq e^{ct}}u(t,x)=1$,

\noindent -  pour tout $c>c_*$, on a $\di\lim_{t\to+\infty}\sup_{\vert x\vert\geq e^{ct}}u(t,x)=0$.
\end{theoreme}
Ce type d'invasion exponentielle a été remarqué par Hamel et Roques \cite{HR} dans le cadre de l'équation (\ref{e2.1}) (i.e. quand $A$ est le laplacien standard) 
pour des données initiales plus grandes que une fonction positive à décroissance algébrique vers $0$ a l'infini.

Notre résultat suivant concerne les solutions avec donnée initiale croissante.
\begin{theoreme} 
\label{t1.2}
Soit $\alpha\in(0,1)$, $N=1$, et considérons la quantité $c_{**}=\di\frac{f'(0)}{2\alpha}$.
Soit $u$ solution de  {\rm (\ref{e1.1})} et supposons que $u(0,\cdot )\not\equiv 0$ est croissante et telle que $\R_-\cap{\rm{supp}}(u(0,\cdot ))$ soit compact. Alors:

\noindent - pour tout $c<c_{**}$, on a $\di\lim_{t\to+\infty}\inf_{x\geq -e^{ct}}u(t,x)=1$,

\noindent -  pour tout $c>c_{**}$, on a $\di\lim_{t\to+\infty}\sup_{x\leq -e^{ct}}u(t,x)=0$.
\end{theoreme}
On montre en particulier que (\ref{e1.1}) ne peut pas admettre d'onde progressive si $\alpha<1$, en contraste avec le 
cas $\alpha=1$. 

On remarquera une différence supplémentaire avec le cas $\alpha=1$: le théorème de KPP pour
$\alpha=1$ assure que la vitesse
d'invasion de l'état stable dans le cas d'une donnée initiale croissante est la même que pour une donnée à support compact. Dans le cas $\alpha<1$, on a deux vitesses différentes.

Le dernier résultat concerne une asymptotique plus précise de la dynamique des lignes de niveau.
\begin{theoreme} 
\label{t1.3}
Soient $\alpha=\di\frac12$, $N=1$ et $f(u)=u(1-u)$.
Soit $u$ solution de  {\rm (\ref{e1.1})} et supposons que $u(0,\cdot )\not\equiv 0$ est à support compact. 
Pour $\lambda\in(0,1)$, soient 
$$
x_\lambda^-(t)=\inf\{x:\ u(t,x)=\lambda\}\ \ \ \hbox{et}\  \ \ x_\lambda^+(t)=\sup\{x:\ u(t,x)=\lambda\}.
$$
Alors, il existe une constante $C_\lambda>1$ tel que, pour tout $t>0$ assez grand,
$$
\di{-C_\lambda e^{t/2}\leq x_\lambda^-(t)\leq -\frac1{C_\lambda}e^{t/2}}\ \ \ \hbox{et}\ \ \ \di{\frac1{C_\lambda}e^{t/2}\leq x_\lambda^+(t)\leq C_\lambda e^{t/2}}.
$$
\end{theoreme}
Ce résultat permet de poursuivre l'étude de la dynamique complète de (\ref{e1.1}). Il admet sans doute une généralisation
pour tout $\alpha\in(0,1)$, que nous ne savons pas effectuer à l'heure actuelle.

\section{Eléments de preuves}
\noindent Dans le cas $\alpha=1$, une façon de faire apparaître la quantité $c_{*,1}$ est la suivante: puisque $f$ est concave,
une sur-solution de (\ref{e2.1}) est la solution $\overline{u}(t,x)$ de
$$
\overline{u}_t-\Delta \overline{u}=f'(0)\overline{u},\ \ \ \ \overline{u}(0,x)=u(0,x).
$$
Examinant le cas particulier $u(0,\cdot )=\delta_0$ - la masse de Dirac en 0 - on obtient 
$
\overline{u}(t,x)=(4\pi t)^{\frac{-N}2}e^{f'(0)t-\frac{\vert x\vert2}{4t}},
$
et $\overline{u}=\lambda$ si $\vert x\vert =2\sqrt{f'(0)}t+o(t)$. Ces considérations très simples permettent aussi de comprendre
pourquoi une décroissance algébrique de la donnée initiale - comme dans \cite{HR} -
conduit aussi à une invasion exponentielle en temps.

Utilisons cette heuristique - déjà effectuée dans \cite{MVV}, par exemple - dans le cas $\alpha<1$. Cette fois-ci, la solution $\overline{u}(t,x)$ de
$$
\overline{u}_t+A\overline{u}=f'(0)\overline{u},\ \ \ \overline{u}(0,\cdot )=\delta_0
$$
est 
$$
\overline{u}= e^{f'(0)t} p(t,x). 
$$
L'estimation
(\ref{e1.2}) indique que $\overline{u}=\lambda$ si 
\begin{equation}
\label{e3.1}
\vert x\vert =O(t^{\frac{1}{N+2\alpha}})e^{c_*t}, \qquad c_*=\di\frac{f'(0)}{N+2\alpha}. 
\end{equation}
La propriété de sur-solution donne une borne supérieure sur la vitesse des
ensembles de niveaux et, partant, la deuxième partie du Théorème \ref{t1.1}. L'obtention  - plus difficile -  de l'autre partie (i.e. la convergence vers~$1$) se fait en deux temps.
D'abord, on montre l'existence de $\varepsilon\in (0,1)$ tel que, si 
$\sigma < f'(0)/(N+2\alpha)$, alors 
\begin{equation}
\label{e3.2}
\liminf_{t\to+\infty}\inf_{\vert x\vert\leq e^{\sigma t}}u(t,x)\geq\varepsilon.
\end{equation}
Ceci est accompli en construisant de solutions de l'équation 
$\overline{v}_t+A\overline{v}=\delta^{-1}f(\delta)\overline{v}$ 
qui prennent des valeurs en $(0,\delta)$ - qui sont en conséquence des sous-solutions de (\ref{e1.1}).
Pour ceci, on résout l'équation pour $\overline{v}$ avec une donnée initiale
prenant des valeurs dans $(0,\varepsilon)$ (avec $\varepsilon < \delta$) jusqu'au temps $T$
où la solution prend la valeur $\delta$, puis en coupant $\overline{v}(T,\cdot)$
par $\varepsilon$, et en
utilisant de façon récursive cet argument.
La convergence vers $1$ est prouvé en utilisant (\ref{e3.2}) et une sous-solution 
liée à la solution de l'équation	
homogène (\ref{e1.0}) avec donnée initiale $\vert x\vert^\gamma$, $\gamma
\in (0,2\alpha)$.

La démonstration du Théorème \ref{t1.3} sur les ensembles de niveaux des solutions
repose sur des arguments similaires combinées avec l'utilisation des
sous et sur-solutions plus précises - dans ce cas $\alpha=1/2$ et $N=1$ - de la forme
$$
\overline{u}= a\left\{ 1+\di\frac{\vert x \vert^2}{b(t)^2} \right\}^{-1},
$$
pour certaines constantes $a$ et fonctions $b=b(t)$ du temps.
Notons que ce Théorème  \ref{t1.3} établit que la position du front 	
donnée par l'argument heuristique  (\ref{e3.1}),
i.e. $\vert x\vert =O(t^{1/2})e^{t/2}$,
n'est pas correct, et à la place est donnée par $\vert x\vert =O(1) e^{t/2}$.

Tous ces points seront développés dans \cite{CR}.



\end{document}